# HISTORICAL ACCOUNT AND ULTRA-SIMPLE PROOFS OF DESCARTES' RULE OF SIGNS, DE GUA, FOURIER AND BUDAN'S RULES

by MICHAEL BENSIMHOUN
26th July 2016, Jerusalem

ABSTRACT

It may seem a funny notion to write about theorems as old and rehashed as Descartes' rule of signs, De Gua's rule or Budan's. Admittedly, these theorems were proved numerous times over the centuries. However, despite the popularity of these results, it seems that no thorough and up-to-date historical account of their proofs has ever been given, nor has an effort been made to reformulate the oldest demonstrations in modern terms. The motivation of this paper is to put these strongly related theorems back in their historical perspective. More importantly, we suggest a way to understand Descartes' original statement, which yet remains somewhat of an enigma. We found that this question is related to a certain way of counting the alternations and permanences of signs of the polynomial coefficients, and may have been the convention used by Descartes. Remarkably, this convention not only provides a ultra-simple proof of Descartes' rule, but it can also be used to simplify the proofs of other related theorems, as the four theorems from the title of the article. Without claiming to be exhaustive, we shall present in this paper an historical account of these theorems and their proofs, and clarify their mutual relation. We will explain how a suitable convention can help understand the original statement of Descartes and greatly simplify its proof, as well as the proofs of the above-mentioned theorems. With the exception of the proof of Fourier's theorem and its generalizations, which run on rudiments of infinitesimal calculus (Taylor's theorem), the proposed demonstrations are so short and elementary that they could be taught at the undergraduate level.

## 1. Historical Perspective

### 1.1. Descartes' rule

Of the theorems listed in the title of this paper, the oldest and by far most famous theorem is Descartes' rule of sign. It was first formulated by Descartes in 1637 in his *Geometry* ([10]). At the beginning of his exposition, the author gives numerical examples of products of polynomials by $X - \alpha$. Next, having claimed that a polynomial $P$ has a root $\alpha$ if and only if it is divisible by $X - \alpha$, he added without proof: "As a result, it is possible to know how many true roots[1] and false roots[2] an equation can have[3]. Namely, it can have as many true roots as the signs + and −−[4] alternate, & as many false roots as two signs + or two signs −− follow one another." [5]

This statement of Descartes' was attacked by several of his contemporaries, who pointed out that a real polynomial can have fewer positive roots than the number of alternations of signs contained in its coefficients. This counter argument is rather surprising, for it is clear from the words of Descartes that he meant a polynomial has *at most* as many positive roots as it contains alternations of signs (actually, this is exactly what Descartes objected in his 77th refutation letter, directed against Roberval ([11])).

On the other hand, other geometers following Descartes attributed this rule to Harriot, an English geometer postdating Viètes and predating Descartes. This attribution is often put to the credit of Wallis, a compatriot of Harriot. Nevertheless, according to Gagneux (citing Stedall), the first who attributed the rule to Harriot seems to be Leibniz, who probably misunderstood the commentary of Wallis on the rule in his treatise of algebra ([16]); it is a fact that this rule does not appear in the work of Harriot, and that it is nowhere revendicated by Wallis for Harriot in his treatise ([16]]). This misunderstanding, as well as the somewhat polemic tone of wallis, caused a sort of franco-english controversy. However, as pointed out by De Gua in his mémoires ([20]), Descartes' rule of signs can very hardly be attributed to Harriot (the arguments of De Gua are well documented and extremely convincing).

Undoubtedly, there is a logical flaw in Descartes' original statement. It seemingly contains two assertions, that should be understood in the following manner:

(1) *A real polynomial has* no more *positive roots than alternations of signs between two consecutive coefficients.*

(2) *A real polynomial has* no more *negative roots than permanences of signs between two consecutive coefficients.*

The following question arises: what is the meaning of the terms *consecutive coefficients*, or in the words of Descartes, "sign that follows one another". Of course, this is clear if the polynomial is not lacunary, but consider for example the polynomial $P(X) = +X^2 - 1$, lacunary in $X^1$. If the term "consecutive" is to be understood in the most obvious manner, then $P$ has exactly one alternation of signs, and no permanence of signs. Therefore, according to assertion (2) above, $P$ should have no negative roots. This is obviously false since $P(X) = (X-1)(X+1)$. Could so trivial a counter-example have escaped Descartes? Of course, one can argue that, unlike Fermat, Descartes was not always fastidious regarding the statement of his theorems; he may, after all, have made a mistake. Nevertheless, he was an excellent geometer and algebraist, whose genius profoundly influenced mathematical thought. Could the mathematician whom even Fermat held in great esteem have missed so simple a point? As we shall show below, this may not be the case. We realized that Descartes' statement is in fact correct if one assigns, in any manner one pleases, signs to the lacunary coefficients. For example, in the aforementioned case, we can write symbolically $P(X) = +X^2 + 0 \cdot X - 1$; we see that $P$ has one alternation and one permanence of signs, which corresponds, indeed, to its unique positive and negative root resp. The same would be true, in this case, if we had written $P(X) = +X^2 - 0 \cdot X - 1$, but generally, the number of alternations and of permanences of signs is not the same if one assigns signs to the lacunary coefficients in different ways. Thus, in order to find the strictest limits of the number of positive and negative roots, and to apply the rule to its greatest extent, it is natural and judicious to attach signs to the lacunary coefficients in such a way that the number of alternations or of permanences of signs is minimized. This can be performed simply by subscribing to a suitable convention. Furthermore, this has the effect of making the proof of Descartes' rule quite simple, and may reasonably be the basis of the proof Descartes had in mind. This convention will be made explicit in the next section, and a ultra-simple proof will be derived. For the moment, however, we continue the historical account.

---

[1] positive roots
[2] negative roots
[3] polynomial
[4] the same as "-"
[5] free translation

As noted above, the meaning of an alternation or a permanence of sign is evident when the polynomial is not lacunary. The situation is more complex when it is. Commonly, lacunary coefficients are ignored during the process of counting the alternations and permanences of signs of $P$. This is the way Descartes' rule was understood by most of the mathematicians from the 19th century up to now, and how it appears in virtually all modern textbooks. Indeed, this method of counting the alternations is sufficient to assert and prove part (1) of Descartes' rule. In fact, we shall see below that this way of counting the alternations of signs coincides with our *minimization convention* above, as far as alternations of signs are concerned; but the two methods do not coincide when applied to the counting of the permanence of signs. This may explain why part (2) of Descartes' rule was generally ignored by mathematicians, likely regarded as a false statement.

It is common to find versions of Descartes' rule that agglomerate the genuine rule with other propositions. The first extension, commonly found in old books, is that if all the roots of a polynomial are real, then the number of positive roots is *exactly* equal to the number of its alternations of signs. Curiously, this proposition is seldom found in modern resources such as Wikipedia, although it could give a useful external criterion to determine the signature of a real quadratic form from its characteristic polynomial, or, in other words, to classify conics. It is difficult to attribute this extension to a specific author, since it was discussed already in Descartes' time as the way some of his contemporaries understood his rule. It certainly appears in the work of De Gua ([19], 1741), which also provide a correct proof of this result.

The second extension consists of the fact that the parities of the number of positive roots and the number of alternation of signs are always equal. This proposition, which several authors have used to prove the rule, is in fact of a very different nature. Indeed, Descartes' rule generalizes to polynomial $f$ defined in any ordered fields, as will be shown below, while the extension above arises from analytic properties of $\mathbb{R}$ (intermediate value and Rolle theorems for instance, or the fact that irreducible polynomials over $\mathbb{R}$ are of a degree no higher than 2). We believe that this *accidental* fact should be deinterlaced from Descarte's rule, at least regarding its proof. We were unable to determine with certainty which author was the first author to give this extension, but we did find it mentioned in the "Notices de l'académie du Gard" of 1809 ([31], p. 204), as one of the theorems proved by Lavernède. Unfortunately, it is difficult to obtain the mémoire of Lavernède, and we were unable to check the proof. This extension is also proved in the 1820 dissertation of Fourier ([15]). More precisely, it was part of Fourier's statement of the Budan-Fourier theorem in the aforementioned dissertation; since the Budan-Fourier theorem is itself a generalization of Descartes' rule, and since Fourier was in possession of these theorems as soon as 1796, we shall attribute this extension to Fourier in this paper. Such an attribution should be confirmed, of course, by further research[6].

Several geometers following Descartes tried to prove his theorem. It is not easy to determine who was the first to give a valid proof; in our opinion, the difficulty in this matter is not in vaguely showing that the rule is correct, but in providing a *rigorous* demonstration. It is also noteworthy that Descartes himself provided a hint on how to prove his rule: Indeed, the words "As a result" just following his remark on the divisibility of a polynomial $P$ by $X - \alpha$ ($\alpha$ root of $P$) suggests strongly that he was in possession of an inductive proof, showing that multiplying a polynomial by $X - \alpha$, with $\alpha > 0$, increases the number of its alternations of signs by at least one unity, while multiplying it by $X + \alpha$ increases the number of its permanences of signs by at least one unity. This proposition, known as Segner's rule, was seemingly not given a correct proof (an even not a correct statement) until the 1828 paper of Gauss, as will be explained in a moment.

---

[6] We were informed that this rule may have been stated by the Chinese mathematician Li Rui in the year 1813.

Roughly speaking, there are two categories of proofs of Descartes' rule: Proofs belonging to the first category are called *algebraic*, and use neither geometry nor infinitesimal calculus, but only algebra and obvious properties of $\mathbb{R}$. The advantage of these proofs is that they can be generalized to ordered fields. Proofs belonging to the second category are called *analytic*, and use geometric curves and extrema to obtain information, or derivatives of the first and higher orders to build inductive arguments. Often, they can be generalized to classes of functions beyond polynomials, which presents an advantage over algebraic proofs.

## 1.2. Demonstrations of Descartes' rule

A demonstration of Descartes' rule was published by Segner in 1728 ([28]). Some authors believe that it is the first rigorous proof of the rule ([3, 30]). Unfortunately, it is difficult to obtain this document, written in latin. Gustaf Eneström gave a short summary in German in Bibliotheca Mathematica ([14]), from which we can see the main ideas:

a) Segner first proves (as did De Gua later) that if an equation $a_0 x^n + \ldots + a_n = 0$ has only real roots, then $a_1/a_0 > a_2/a_1 > \ldots > a_n/a_{n-1}$.

b) He then multiplies the equation by $x - m$ and looks at the sign changes in the cases where $m > a_1/a_0$, $a_1/a_0 > m > a_2/a_1$, $\ldots a_n/a_{n-1} > m$.

c) He shows that the multiplication always introduces a new sign change, and that this still holds when the coefficients are not all positive.

d) An equation has thus as many positive roots as the number of sign changes in its coefficients, and the same method can be applied for the negative roots.

Enerström qualifies this proof as "essentially correct, although neither particularly clear nor beautiful", and according to this report, it is restricted to polynomials without imaginary roots.

Segner published another proof of the rule in 1756 ([29]), apparently not restricted to polynomials without imaginary roots. As we shall explain in a moment, the second proof of Segner fails to be correct, unless it be properly corrected and completed. The last proof of Segner was also explained and extended in the Encyclopaedia Britannica in 1824 ([13]), apparently without compensating for the important failure in the proof. According to this entry, Segner's proof would be "not only the most simple, but probably the most simple that will ever be invented." In our opinion, Segner's proof, may be convincing (after correction), but it is certainly not the simplest possible proof produced since.

It is generally acknowledged that the first rigorous demonstration of Descartes' rule was published by De Gua in 1742 ([19]). In fact, De Gua stated two propositions, of which he gave separated demonstrations. The first proposition stated by De Gua is that, if all the roots of a polynomial $P$ are real, the number of positive roots of $P$ is equal to its number of alternations of signs, and the number of negative roots of $P$ is equal to its number of permanences of signs. This is Descartes' rule in the case where all the roots are real. In order to prove this result, De Gua, following the lines sketched by Descartes, proposed to show the following, namely,

*If a polynomial, whose all roots are real, is multiplied by $X - \alpha$ with $\alpha > 0$, then its number of alternations of signs increases by one unity. If, in contrast, it is multiplied by $X + \alpha$, then its number of permanences of signs increases by one unity.*

This is almost the same proposition Segner attempted to show in 1756. From this time, it has been known as Segner's lemma:

*If a polynomial is multiplied by $X - \alpha$, with $\alpha > 0$, then its number of alternations of signs increases by* at least *one unity. If, in contrast, it is multiplied by $X + \alpha$, then its number of permanences of signs increases by* at least *one unity.*

To prove his theorem, De Gua based his argument on a remarkable lemma:

*If a polynomial contains only real roots, then its coefficients $a_i$ fulfill*
$$a_i^2 - a_{i-1} a_{i+1} > 0.$$

De Gua seems to have been unaware that this proposition, albeit presented in a different way, had already been proved by Segner in 1728. As for Segner, he found an ingenious and simple argument based on "passages" from the coefficients of one polynomial to another. But the most salient fact regarding the proof of De Gua and the second proof of Segner (few years later) is the following: both geometers did not pay much attention to lacunary coefficients in polynomials; De Gua, explaining how the alternations and permanences of signs of a polynomial have to be counted, just indicated that lacunary coefficients may be replaced by infinitesimal quantities of arbitrary signs ([19], p 83). Similarly, Segner, in the aforementioned paper, contented himself with saying that lacunary coefficients have to be replaced by $\pm 0$, and added nothing more on this question in the remaining part of his discussion ([29], p. 292). Under this method of counting the alternations and permanences of signs, Descartes' rule of signs holds, as contended above. What does not, in general, is Segner's lemma: Indeed, let us consider, for instance, the polynomial $P(X) = X^3 - 1$. This polynomial, multiplied by $X - 1$, is equal to $Q(X) = X^4 - X^3 - X + 1$. According to the convention of De Gua and Segner, they could be written $X^3 - 0X^2 + 0X - 1$, and $X^4 - X^3 - 0X^2 - X + 1$. Therefore, the number of alternations of $P(X)$ would be 3, while the number of alternations of $Q$ would be 2, in contradiction with Segner's lemma. This simple example shows that not only was the proof of Segner non rigorous, but also that his statement of the lemma was incorrect! In the case of De Gua, nevertheless, things are different, because De Gua stated and proved the lemma under the additional assumption that all the roots of $P$ are real. This implies, in particular, that the lacunary coefficients of $P$ never appear consecutively (as can be seen, for example, from De Gua's lemma above). Thanks to this additional condition, Segner's lemma proves be true, even if one assigns arbitrary signs to lacunary coefficients, as shown by De Gua.

To close this discussion, we point out that a correct statement of Segner's lemma can be given with the help of a suitable convention for counting the signs of the lacunary coefficients inside polynomials, as done further in this article. Once such a convention is assumed, the proof of Segner proves to be perfectly valid, just as other results related to sign alternations, invented in the past.

The second proof of De Gua is analytic, and no assumption is made about the reality of the roots of the equation. It is an application of the geometry of curves to algebra, in the spirit of certain recent proofs (see e.g. [32]). This proof includes three steps: first, the statement that a polynomial $P$, whose coefficients share the same sign, has no positive roots (theorem 1). Second, De Gua observes that if two coefficients $a_k$ and $a_{k-1}$ of $P$ are of opposite signs, then the differential of $F(X) = X^{-k}P(X)$, multiplied by $X^{k+1}$, has exactly one alternation less than $P$ (theorem 2[7]). Third, De Gua shows that if $P(X)$ has $n$ positive roots, counted with multiplicities, then the polynomial $G_1(X) = X^{k+1}F'(X)$, has at least $n-1$ positive roots (theorem 3). With respect to modern standards, the demonstration of this assertion, based on geometric figures and a on vague argument involving "infinitesimals", may seem somewhat insufficient; but it was considered rigorous by the time of De Gua. In fact, using modern tools, the proof is quite simple: It is easy to see that every root $\zeta_i$ of multiplicity $m_i$ of $P$ is a root of multiplicity $m_i - 1$ of $G_1(X)$. On the other hand, between two roots of $P$, which are also zeros of $F$, $F'$ must vanish at least one time according to Rolle's theorem. Therefore the number of positive roots of $G_i(X)$ is at least

$$\sum_i (m_i - 1) + n - 1 = \sum_i m_i - 1 = n - 1.$$

The final step in the proof of De Gua is to use the second and third theorems as many times as necessary, until a polynomial with coefficients of like signs is obtained: Starting from $P$, we saw that $G_1$ has exactly one alternation less that $P$, and at most one root less than $P$. From $G_1$, a polynomial $G_2$ can be derived in the same manner that $G_1$ was derived from $P$; it contains exactly one alternation, and at most one root, less than $G_1$. Continuing this process, one arrives to $G_k$, which has no alternations of signs, exactly $k$ alternations, and at most $k$ roots, less than $P$. But $G_k$ has no positive roots according to the first theorem of De Gua, hence, the number of roots of $P$ is at most equal to the number of alternations of $P$.

---

[7] Actually, this theorem was stated in term of arithmetic progressions, and is essentially equivalent to this assertion.

It is surprising that this ingenious proof of De Gua was not simplified and reshaped into a short inductive form by mathematicians of the first half of the 19th century. Even Gauss felt the need to invent another proof, and did not try to retranscribe that of De Gua in modern terms. In fact, De Gua's method of proving Descartes' rule was rediscovered by Laguerre in 1883 ([5]), and more recently by Komornik in 2006. The extremely simple proof of Laguerre was considered to be entirely new by his contemporary (see the preface of Poincaré in [5]), despite it is essentially the same as De Gua's proof. Also, Komornik was seemingly unaware of the proof of Laguerre, since it published exactly the same proof ([23]). But De Gua's contribution to algebra is not limited to these proofs. In another *mémoire* ([20]), he gives an extended and precise historical analysis of the theory of equations up to his time, and investigates means to evaluate the number of complex roots of real polynomials. The so-called *De Gua's rule* was extracted from this work, a kind of lower-bound on the number of complex roots of real polynomials (see Sec. 2.3 below). In conclusion, we believe that it is not exaggerated to say that the work of De Gua is one of the most valuable piece in the theory of equations, created in the 18th century.

In the years following De Gua and Segner, several authors presented other demonstrations. For example, a laborious and long analytic proof was given by Æpinus in 1758 ([27]), and an algebraic one was given by Lagrange in 1808 ([25]), following the lines of his theory of equations. On the same note, Lagrange mentions a proof of Kæstner, which we were unable to find. Regarding the proof of Æpinus, it is based on several propositions involving differentials of polynomials of orders 1 or more, and is somewhat difficult to follow. The author does not enunciate what he calls "Harriot theorem", and in particular, if he assumes that all the roots are real. Even more troublesome, he does not explain how he envisages to handle the alternations and permanences of signs in lacunary polynomials. As a result, it appears that the proof of Æpinus cannot be considered as rigorous (or even correct). But the argument of Æpinus to prove Descartes' rule, based on differentials of any order, predated the method of Fourier regarding this kind of questions.

In 1828, an extremely simple demonstration of Descartes' rule was published by Gauss ([17]). It is exposed in a somewhat more understandable form in [12]. May Gauss have been motivated by an article of Grunert one year earlier ([18]), in which this later criticized the previous demonstrations of the rule, and proposed his own lengthy and complicated proof? Despite its simplicity, the proof of Gauss does not seem to be well known, even by renowned mathematicians (in [12], Dickson said the proof was communicated to him by Curtiss, and seemed to ignore it is due to Gauss). This may explain why proofs for Descartes' rule are still being invented. Curiously, Abraham Adrian Albert was unsatisfied with the proof of Gauss, and produced a far more complicated proof in 1943 ([1]). The neatly simple proof of Gauss is algebraic, following, once more, Descartes' strategy above, and is based on sign alternation considerations only. But the most important thing is that, for the first time (apparently), the function of lacunary coefficients in Descartes' rule was fully understood and exhibited: First of all, Gauss explains how to count the alternations of signs; he simply proposes to "ignore" lacunary coefficients, and to order polynomials by descending powers of $X$. According to Gauss, an alternation of sign is simply a sign change between two consecutive powers of $X$ (this simple and natural way of counting alternations of signs has become the standard). Armed with this convention, Gauss was able to give a correct formulation of Segner's lemma, and to produce an extremely simple and elegant proof. In fact, his demonstration proves to be essentially the same as that of Segner's, but correctly stated, and formulated with the illuminating clarity of pure algebra. Regarding the permanences of signs of $P(X)$, Gauss simply regards them as alternations of $P(-X)$. Thus, he can provide a direct definition of sign permanences. This *ad hoc* method has the merit to overcome the problem of lacunary coefficients inside polynomials. Nevertheless, it is possible to give a less artificial (but ultimately equivalent) convention, as set out further in this paper.

In the same article, Gauss explicitly states the first of the two extensions mentioned above, together with an assertion involving (implicitly) De Gua's rule:

*If one counts the total number of missing powers inside $P$, and then for each odd length gap, decreases it by one if it separates a sign change, and increases it by one if it separates a sign permanence, one thereby obtains a number which the number of imaginary roots of $P$ must at least equal.*

This rule is in fact equivalent to Cor. 2.3.1 below, a simple corollary of Descartes' rule.

In 1883, an important step was accomplished: Laguerre published an ultra-simple, seemingly new, analytic proof of Descartes' rule ([5], p. 3, [23]). As pointed out above, it is in fact essentially the same as De Gua's proof, but by the time of Laguerre, it was possible to formulate it in extremely clear and concise terms. For example, the notion of multiple root needs not be restricted to polynomials: it extends immediately to any analytic function $f$, as being the order of the first differential of $f$ that does not vanish at a zero of $f$. Also, Gauss's method for counting the alternations of signs of polynomials had become standard. The proof of Laguerre (or rather his formulation) can be explained in few words: Let us call *extended polynomial* a function $P$ defined in $\mathbb{R}_+$, of the form
$$P(X) = a_{\lambda_1} X^{\lambda_1} + a_{\lambda_2} X^{\lambda_2} + \cdots + a_{\lambda_p} X^{\lambda_p},$$
where $\lambda_i \in \mathbb{R}$, $\lambda_i > \lambda_{i+1}$ and $a_{\lambda_i} \neq 0$. The idea of Laguerre's is to prove that Descartes' rule is in fact valid for extended polynomials. Of course, this implies its validity for ordinary polynomials. It is a typical example of the fact that it is sometimes simpler to prove an extension of a theorem than the theorem itself. Indeed, this assertion can be assumed inductively for extended polynomials with $m$ alternations of signs (it is evident if $m = 0$). Let
$$P(X) = a_{\lambda_1} X^{\lambda_1} + a_{\lambda_2} X^{\lambda_2} + \cdots \quad (a_{\lambda_i} \neq 0)$$
be an extended polynomial, with $n$ roots and $m+1$ alternations of signs. For some index $i$, holds $a_{\lambda_i} a_{\lambda_{i+1}} < 0$. Laguerre applies De Gua's trick: multiplying $P$ by $X^{-\lambda_i}$ leads to an extended polynomial
$$Q(X) = a_{\lambda_1} X^{\lambda_1 - \lambda_i} + \cdots + a_{\lambda_i} X^0 + a_{\lambda_{i-1}} X^{\lambda_{i-1} - \lambda_i} + \cdots.$$
It is clear that the number of alternations of signs of $Q$ is the same as that of $P$, as are the numbers of roots of $P$ and $Q$ in $\mathbb{R}_+$. Moreover, it is easily seen that the differential $Q'$ of $Q$ has exactly one alternation of signs less than $P$. Therefore, to complete the induction, it suffices to show that the number of roots of $Q'$ in $\mathbb{R}_+$ is at least equal to $n - 1$ (counting roots with their multiplicities). The 3 lines proof of this assertion can be seen in the discussion about De Gua's proof above.

What is even more impressive is that the proof of Laguerre's generalizes immediately to infinite series of the form
$$\sum_i a_{\lambda_i} X^{\lambda_i},$$
with a finite number of alternations of signs. Thanks to his extension of Descartes' rule, Laguerre was able to prove another powerful theorem of his own, which seems to have fallen into oblivion ([5], p. 6):

*Let the following sequence of polynomials be considered*
$$\begin{aligned} P_0(X) &= a_n X^n + a_{n-1} X^{n-1} + \cdots + a_0 = P(X), \\ P_1(X) &= a_n X^{n-1} + a_{n-1} X^{n-2} + \cdots + a_1, \\ &\cdots \\ P_n(X) &= a_n, \end{aligned}$$

*and let $a$ be a positive number. Then the number of roots of $P$ larger than $a$ is equal to the number of alternations of signs in the sequence*
$$P_0(a), P_1(a), \ldots, P_n(a),$$
*or is larger than it by an even number.*

We point out that if $a = 0$, it is nothing more than Descartes' rule. But this theorem enables to find an upper bound on the roots of $P$ larger than $a$ more easily than transforming $P$ into $P(X - a)$, and applying subsequently Descartes' rule. Indeed, $P_i(a)$ can be computed inductively:

$$P_i(a) = aP_{i-1}(a) + a_i.$$

In the aforementioned article, Laguerre was even able to give a theorem parallel to Budan's rule, in the spirit of the theorem above. ([5], p. 9–11). Further, Laguerre gave numerous applications of these results and other theorems, beyond the scope of this paper. In another article ([5], p. 48-50), Laguerre proved several results related to homogenized polynomials, among which an extension of Rolle's theorem to complex roots. Ending his paper, he asked for a generalization of Descartes' rule of signs that could handle complex roots. This question, and other ones by Laguerre, widely stimulated mathematicians which followed Laguerre (see e.g. [7, 8]), in particular Nicholas Obrechkoff who used Cauchy's theorem to obtain an extension to roots inside a cone ([26]).

A ultra-simple algebraic proof of Descartes' rule was finally published by Krishnaiah in 1963 ([24]), based, again, upon Descartes' strategy (Segner's lemma). In his paper, Krishnaiah calls the proof of Segner's lemma which appears in the literature "a diagrammatic persuasion, whose rigorous presentation is rather lengthy", and proposes a rigorous inductive proof to compensate for this deficiency. His opinion is, at the least, difficult to understand: One does not see in what the proof of Gauss or that of Laguerre fails to be rigorous. The guiding principle of Krishnaiah's proof, which is in essence the same as Gauss's, is the following: assume inductively that multiplying a polynomial with $r$ alternations of signs ($r > 0$) increases its number of alternations by one unity (Segner's lemma). Consider a polynomial $P = a_0 + \cdots + a_n X^n$ containing $r+1$ alternations of signs, and let $a_s$ be the last coefficient at which an alternation of sign occurs (say $a_s a_{s+1}$). Write the polynomial in the form $Q + R$, where $Q$ is the part of $P$ up to the $s$th coefficient, and $R$ is the remaining part of $P$. In particular $Q$ contains $r$ alternations of signs. Now, $(X - \alpha)P = (X - \alpha)Q + (X - \alpha)R$. It is important to note that because of the definition of $s$, the last and leading coefficients of $R$ share the same signs. Hence it is not difficult to show that $P$ must contain at least $r + 2$ alternations of signs, completing the induction.

Since the proofs of Gauss and Krishnaiah, as well as the tentative proof of Descartes presented in a further section of this paper (Sec. 2.1), are probably the simplest algebraic proofs that will ever be invented, it may be interesting to compare them. Gauss's proof is non-inductive, and may be preferred over the other proofs by persons experienced in algebra, for it allows to embrace the intrinsic truth of Descartes' rule in a glance. The proof of Krishnaiah can be seen as a formulation of the proof of Gauss in an inductive form; this makes it simpler for persons not experienced in algebra, and also somewhat more rigorous, though it may be considered inferior because of its inductive nature. Both proofs do not use any special conventions to count the alternations of signs of polynomials: lacunary coefficients are simply ignored. Interestingly, in the proof of Krishnaiah, the induction is performed on the number of alternations of signs, and not on the degree of $P$. Not so in the tentative proof of Descartes below: There, the induction is performed on the degree of $P$, which is perhaps more natural. As mentioned above, what prevents most of the proofs in the old literature from being rigorous is the lack of a suitable convention for counting the number of alternations and permanences of signs of the polynomial coefficients. The merit of the proofs of Gauss and Krishnaiah is that they overcome this difficulty without introducing any special convention or external concepts. Nonetheless, once a suitable convention is established, Descartes' rule follows in a straightforward manner (see Sec. 2.1 below), as does the extensions outlined in the next section: De Gua's rule can be shown to be an artifact of this convention, while what is referred to below as "Fourier's rule" appears to be a simple corollary of complex analysis, deinterlaced from Descartes' and Budan's rules. Moreover, even the proof of Budan's rule is simplified by this convention.

In conclusion, the proof of Krishnaiah may be preferred if the goal is to provide a quick and rigorous proof of Descartes' rule, free from any other considerations. The tentative proof of Descartes below may be preferred over that of Gauss for its historical perspective and its suitability to be taught at a very elementary level. It may also be preferred over Krishnaiah's proof if a deeper insight of the principles involved is desired, in order to derive easily extensions of Descartes' rule.

And what about the analytic proof of Laguerre? Apart from being extremely simple, it immediately provides an extension of Descartes' rule to polynomials with real coefficients, as well as to infinite series with a finite number of alternations of signs. Even though it cannot be taught at a very elementary level, it is certainly the best choice for persons having some experience in calculus.

### 1.3. Budan and Fourier's rules

Between the years 1796 and 1822, Descartes' rule took an unexpected course. Both Budan and Fourier found that it can be generalized in such a way that it provides an upper bound of the number of positive roots of a polynomial, *between* any two given bounds $a$ and $b$. Namely, if one counts the number of alternations of signs inside the sequences $P(a), P'(a), \ldots, P^{(n)}(a)$ and $P(b), P'(b), \ldots, P^{(n)}(b)$, then the first number is always in excess over the second one by a quantity equal to the number of roots of $P$ inside $]a,b]$, or exceeding it by an even number. In the literature, this theorem is referred to as Budan's rule, Budan's theorem, Fourier's theorem, or the Budan-Fourier theorem. Actually, Budan formulated this theorem in terms of polynomial coefficients only, while Fourier used the differentials of $P$. These two formulations are, of course, equivalent, according to Taylor's Theorem. It is evident that both Budan and Fourier knew Taylor's theorem, and hence, that this theorem could be expressed in both forms. But Fourier was more concerned with the theoretical aspect of this question, and with the generalization to other classes of functions, while Budan was more concerned with the *algorithmic* aspect, and in particular with the root separation problem for polynomials. It is interesting that the method of Budan is now used inside the fastest algorithm for root separation to date, while the formulation of Fourier has inspired a number of theoretical studies.

The credit of discovery of Budan's rule occasioned a dispute of sorts between the partisans of Fourier and those of Budan. Darboux gave an historical account of the work of the two authors in a note ([9]). Since Budan, Fourier and Darboux were all French, the latter cannot be suspected of partiality. The salient facts concerning the paternity of the Budan-Fourier theorem will be now exposed.

During the years 1796-1803, Fourier taught the rule at the *École polytechnique*, and was most overlooked at the time of Budan. Nevertheless, he did not published anything, and several communications to the *Institut de France* were lost. Budan stated his rule in 1806 in a treatise, but was unable to prove it. In 1811, he presented a *mémoire* on the resolution of equations, including an algebraic proof of the rule, to the *Académie des Sciences*. This proof was found to be correct by the commissars Lagrange and Legendre. Actually, since it is based on Segner's lemma, it is affected by the same deficiency as the lemma[8], but this deficiency seems to have been unnoted by the time of Budan. In 1820, Fourier published his paper [15], in which he investigated the theorem and its application to the theory of equations to a greater depth than Budan had done. The analytic principles in the proof of Fourier turned out to be far superior to that of Budan. We point out that neither Fourier, nor Budan, stated explicitly how to handle lacunary terms during the process of counting the alternations and permanences of signs in sequences of numbers. In the case of Fourier, nevertheless, this may be understood from his demonstration, and corresponds essentially to the minimization convention presented in this paper (Sec. 2). In 1822, Budan published a new *mémoire* on the resolution of equations, in which he included the demonstration of the rule he gave in his first 1811 *mémoire* ([4]). In the opinion of Darboux, the rule should be, without doubt, attributed to Fourier. The arguments of Darboux may be convincing, but we believe that the only *objective* criterion to attributing a theorem is that its author has published it, or in the very least, has showed it in a sufficiently wide, public forum. Actually,

---

[8] This was discussed above.

Navier, who published posthumously an unachieved book of Fourier in 1830 ([9], p. ii–iii), asserts, in the preface of the book, that the copy of a treatise of Fourier by Bonnard, is in his hands, and was shown to Roux in 1795. In the last pages of this copy, Roux attests on his honor to have seen the said copy by the year 1795, and that Bonnard said to him that it was written by Fourier when he was only 18 years old. Among the results found in this treatise and described in the preface, Navier quotes the aforementioned rule and its proof. According to this source, Fourier was in possession of the theorem as soon as 1784. Nevertheless, we have to point out the following point: the demonstration of Fourier, like the one published in 1820 or that published posthumously by Navier in 1830, is far from being rigorous. It is true that it contains all the notions needed to prove the theorem, but it is more a laborious persuasion of the author than a proof in the formal meaning of the term. This may explain why Fourier delayed the publication of his investigations on this subject until the years 1800-1820, despite the repeated requests of his colleagues between the years 1800–1820. Fourier had an extraordinary creative genius, but unlike Gauss, he was not good at proving things rigorously. He simply felt the right arguments and tried to use them to persuade the reader. As a result, he may have felt that his arguments were not presented with sufficient rigor[9]. In contrast, the proof given by Budan in 1811 is rigorous enough (at the exception of the use of Segner's lemma, as explained above). Thus, even if Budan's rule was discovered and proved by Fourier for the first time, as seems very probable, we believe it should be formally attributed to Budan. Another alternative, followed by several mathematicians like Akritas, is to make a distinction between "Budan's rule" and "Fourier's theorem". Following this approach, Budan's rule would be the direct extension of Descartes' rule for polynomials, applied between two bounds $a$ and $b$, and expressed in terms of *polynomial coefficients*. On the other hand, Fourier's theorem would be the generalization of Descartes' rule for polynomials (and even for differentiable functions), applied between two bounds $a$ and $b$, and expressed with the help of the *differentials of f at a and b*. This is probably the most precise approach from the historical point of view, and is also fair with respect to the two inventors of this theorem.

It is to be observed that Descartes' rule can be deduced easily from the Budan's, and that it is not only valid for polynomials, but also for functions that admit an $n$-th differential of constant sign inside some interval $]a, b]$ (see Thm. 2.4.1 below). Also, in an interesting paper dealing with recent discoveries in this domain, it was noted by Curtis in 1918 that Budan's rule gives the exact number of roots between $a$ and $b$, if and only if all the roots of the polynomial are real ([8]). A somewhat vague proof of Budan's rule was given in Dickson ([12]). Often, proofs of Budan's rule, and in particular the proofs of Budan's and Fourier's, fail to be sufficiently rigorous to comply with modern standards, even though they would be valid if suitably completed. In particular, we believe that an argument equivalent to the least-upper bound property cannot be avoided in general analytic proofs, valid for functions of the above-mentioned type. Nevertheless, in the case of polynomials, algebraic proofs can also be given, using Descartes' rule or Segner's lemma (this is in fact the way Budan proved his theorem). Another modern and simple algebraic treatment, based on Descartes' rule, can be found in [6]. This last proof is probably the simplest that will ever been invented and, being algebraic, has moreover the merit of being valid for polynomials with coefficients in real fields. However, it cannot be generalized to classes of functions beyond polynomials. In this paper, a rigorous treatment of Fourier's rule will be presented, valid for a more general class of functions.

---

[9] Thm. 2.4.1 below can be seen as a rigorous exposition of the ideas of Fourier. It is based on the compactness principle, which was unavailable by the time of Fourier.

An ultimate and efficient algorithm providing the exact number of positive roots between two bounds, was finally found by Sturm. This famous theorem, presented to the academy of Paris in 1829, was strongly inspired by the work of Fourier, as acknowledged by its author. Sturm's algorithm efficiently solves the problem that has preoccupied mathematicians since Descartes: locating the roots of a polynomial, in order to compute them using known approximation algorithms. An excellent account of the origin and influence of Sturm's algorithm in mathematics, as well as its deep implications as found by Tarski, is given in [21].

During decennia, Sturm's algorithm completely eclipsed the Budan-Fourier theorem. But in 1911, Hurwitz published a proof of a beautiful extension of this theorem ([22]).

*Let $f$ be a real valued and infinitely differentiable function, defined on an interval $[a, b]$. Suppose that for some $m \geq 0$, $f^{(m)}(x) \neq 0$ at $x = a$ and $x = b$. Let $V(x, m)$ denotes the number of alternations of signs in the sequence*

$$f^{(m)}(x), f^{(m-1)}(x), \ldots, f'(x), f(x),$$

*and $Z(f)$ and $Z(f^{(m)})$ denote the number of zeros of $f$ and $f^{(m)}$ in $]a, b]$ resp. Then $V(a, m) - V(b, m)$ is equal to $Z(f) - Z(f^{(m)})$, or exceeds it by an even positive number.*

This result of Hurwitz, that has almost fallen in oblivion, can be extended to a more general class of functions than those considered by Hurwitz. This is the object of Thm. 2.5.2, presented at the end of this paper.

Although Sturm's algorithm completely eclipsed it for decennia, the Budan-Fourier theorem has seen a renewal of interest in the last decades, since it lies at the heart of (through the Vincent algorithm) the most powerful method of polynomial root isolation: the VAS algorithm (2005).

## 2. Descartes' Presumed Convention and its Applications

This section is devoted to the exposition of the convention mentioned above, and to showing how it leads to a very simple and natural proofs of Descartes' rule. It will also be shown how it simplifies the proof of other results such as De Gua's rule and the Budan-Fourier theorem. From an historical point of view, this may reflect Descartes' thought more precisely than other methods. The proof of Descartes' rule, given in this section, has some similarities with [2], but it is simpler and free of unnecessary analytic considerations. As a consequence of being purely algebraic, it remains valid for polynomials with coefficients in real field.

For the sake of completeness, Fourier's rule and Hurwitz's extension of Fourier's theorem will also be presented at the end of this section.

CONVENTION: Consider a polynomial $P(X) = a_n X^n + a_{n-1} X^{n-1} + \cdots + a_0$, where some coefficients $a_i$ may be null. If two contiguous coefficients $a_{i+1}$ and $a_i$ are both positive or both negative, then the pair $(a_{i+1}, a_i)$ contributes 1 permanence and 0 alternation of signs. If one of them is positive and the other is negative, then it contributes 1 alternation and 0 permanence. If now the polynomial has one or more coefficients equal to 0, then in the context of counting the alternations of signs of $P$, every coefficient $a_i = 0$ is considered to be of the same sign as $a_{i+1}$. Therefore, the sign of a non-zero coefficient limiting a sequence of null coefficients on the left propagates to the entire sequence. For example, if $P(X) = 3X^4 - X$, we assign the following signs to the lacunary coefficients: $P(X) = 3X^4 + 0 \cdot X^3 + 0 \cdot X^2 - X - 0$. In practicality, this is the same as ignoring lacunary coefficients.

On the contrary, in the context of counting the permanences of signs of P(X), we consider the sign of a coefficient $a_i = 0$ to be opposite to the sign of $a_{i+1}$. Thus, the signs of a sequence of null coefficients alternate, starting from the first non-zero coefficient limiting the sequence on the left. For example, using the previous polynomial, the signs are $P(X) = 3X^4 - 0 \cdot X^3 + 0 \cdot X^2 - X + 0$. It is important to note that unlike alternations of signs, this convention does not amount to dropping the null coefficients. In fact, the reader should mentally check that the following rule holds: a sequence of coefficients of the form $a_i, 0, 0, \ldots, 0, a_j$, with $a_i a_j \neq 0$, contributes 1 permanence if $a_i$ and $a_j$ are of opposite signs and the number of "0" is odd; or, if $a_i$ and $a_j$ are of the same signs and the number of "0" is even. Otherwise, it contributes no permanence.

It may seem strange to use two different methods in order to count the number of alternations and of permanences of signs, but the advantages this method offers are twofold: First, with this convention, the number of alternations and sign permanence is always minimized in lacunary polynomial, as can be easily verified. Second, these two methods are dual in the sense that the number of permanences of signs of $P(X)$ is always equal to the number of alternations of signs of $P(-X)$. The verification of this simple fact is left to the reader, as well.

In the rest of this section we shall show how this convention can be applied in order to provide a ultra-simple proof of Descartes' rule, and to simplify the demonstrations of the other rules cited in the title of this paper.

NOTATIONS: Henceforth, we use $z^+(P)$ to denote the number of positive roots of a polynomial $P$, $z^-(P)$ to denote the number of negative roots, and $z^0(P)$ to denote the number of null roots of $P$; roots are always counted with their multiplicities. Of course, the number $z^0(P)$ can be *read over* from $P$ itself, as it is the largest power of $X$ dividing $P$. We also use $v(P)$ to denote the number of alternations of signs of $P$, and $c(P)$ the number of its permanences (hence $c(P(X)) = v(P(-X))$). We say that a coefficient $a_i$ is *trailing* if $a_i = a_{i-1} = \cdots = a_0 = 0$.

We note that if the polynomial $P$ is not lacunary, except, perhaps, in its trailing coefficients, then
$$v(P) + c(P) = \deg(P) - z^0(P). \tag{1}$$

Indeed a pair of two successive coefficients is either in alternation or in permanence of signs. On the other hand, if $P$ is lacunary, then according to our convention, outlined above, a block of successive coefficients of the form $a_i, 0, 0, \ldots, 0, a_j$ (with $a_i a_j \neq 0$) contributes at the most one alternation and one permanence. Since any such block contains at least 3 coefficients (and therefore at least two pairs of successive coefficients), the sum of the numbers of alternations and permanences of $P$ is not larger than the overall number of pairs of successive (and not trailing) coefficients in $P$. In other words, it is not larger than $\deg(P) - z^0(P)$:
$$v(P) + c(P) \leq \deg(P) - z^0(P). \tag{2}$$

## 2.1. Descartes' Rule of Signs

With the notations above, Descartes' rule and its special case can be formulated as follows:

*For any real polynomial $P$,*
$$v(P) \geq z^+(P) \quad and \quad c(P) \geq z^-(P).$$

*Furthermore, if all the roots of $P$ are real, then*
$$v(P) = z^+(P) \quad and \quad c(P) = z^-(P).$$

*Proof*: The first assertion of the theorem implies the second, since if all the roots of $P$ were real and such that $v(P) > z^+(P)$ or $c(P) > z^-(P)$, there would hold
$$\deg(P) - z^0(P) = z^+(P) + z^-(P) < v(P) + c(P),$$

contradicting $(2)_{\text{p.17}}$.

To prove that $v(P) \geq z^+(P)$, note that if $\alpha$ is a real root of $P$, then

$$P = (X - \alpha)Q(X)$$

where $Q$ is a polynomial of degree smaller than $P$. Thus, by an evident induction on the degree of $P$, it suffices to prove the so called Segner's lemma (Descarte's strategy):

*Multiplying a polynomial $Q(X)$ by $X - \alpha$, with $\alpha > 0$, increases the number of its alternations of signs.*

We can see that $c(P) \geq z^-(P)$ follows from $v(P) \geq z^+(P)$ by changing $P(X)$ into $P(-X)$. Let

$$Q(X) = a_n X^n + \cdots + a_k X^k \quad \text{and} \quad P(X) = (X-\alpha)Q(X) = b_{n+1}X^{n+1} + \cdots + b_k X^k,$$

where $k \geq 0$ and $a_k \neq 0$. The following relations hold:

$$b_i = \begin{cases} a_{i-1}, & \text{if } i = n+1, \\ a_{i-1} - \alpha a_i, & \text{if } k < i \leq n, \\ -\alpha a_i, & \text{if } i = k. \end{cases} \quad (3)$$

Let us denote the sign of $a_i$ by $s_i$, and the sign of $b_i$ by $s'_i$:

$$s_i, s'_i \in \{+, -\}.$$

We input the signs of $a_i$ and $b_i$ into a "table of signs" in the following manner:

| i | $n+1$ | $n$ | $n-1$ | $\cdots$ | $k$ |
|---|---|---|---|---|---|
| sign of $a_i$ |  | $s_n$ | $s_{n-1}$ | $\cdots$ | $s_k$ |
| sign of $b_i$ | $s'_{n+1}$ | $s'_n$ | $s'_{n-1}$ | $\cdots$ | $s'_k$ |

For the sake of simplicity, let us denote such a table by

$$(s_n, s_{n-1}, \ldots, s_k; s'_{n+1}, s'_n, \ldots, s'_k),$$

and call the number $n - k$ the *extent* of the table.

Because of relation (3), any such table should satisfy the following three properties:

1. $s'_{n+1} = s_n$
2. $s'_k \neq s_k$,
3. if $s_i \neq s_{i-1}$, then $s'_i = s_{i-1}$.

It is important to note that these properties hold even if $P$ and $Q$ have lacunary coefficients, due to the convention above. To prove Descartes' rule, it suffices to prove that the third row in the table contains more alternations of signs than the second row. To this end, we need nothing more than the three aforementioned properties.

If the extent of the table is $1$, this is obviously true by the first and second properties. Let us inductively assume this assertion for tables of extent of at most $n - k - 1$. Let $T = (s_n, s_{n-1}, \ldots, s_k; s'_{n+1}, s'_n, \ldots, s'_k)$ be a table of extent $n - k$.

Assume first that no alternation of signs occurs in the sequence $(s_n, \ldots, s_k)$ (that is, it is constant). Then $s'_{n+1}, \ldots, s'_k$ must contain at least one alternation of signs: indeed, by hypothesis, $s'_{n+1} = s_n$ and $s'_k \neq s_k$, hence $s'_{n+1} \neq s'_k$. Thus, there must exist some $i \in \{1, \ldots, n+1\}$ such that $s'_i \neq s'_{i-1}$. This shows that in this case, the induction hypothesis is fulfilled for table $T$.

Now assume that an alternation of signs occurs at some rank $i > k$ inside the sequence $s_n, \ldots, s_k$:

$$s_i \neq s_{i-1}.$$

By the third property above, $s'_i = s_{i-1}$, hence $s'_i \neq s_i$. We can extract from $T$ the following table, of smaller extent:

$$T' = (s_n, \ldots, s_i; s'_{n+1}, \ldots, s'_i).$$

This table obviously satisfies the three properties above. Therefore, using $A_1$ to denote the number of alternations of the sequence $s_n, \ldots, s_i$, and $A'_1$ to denote that of $s'_{n+1}, \ldots, s'_i$, the induction hypothesis implies $A'_1 > A_1$. On the other hand, we can extract another table of smaller extent from $T$, namely

$$T'' = (s_{i-1}, \ldots, s_k; s'_i, \ldots, s'_k).$$

It satisfies the properties above as well, hence the sequence $s'_i, \ldots, s'_k$ must also contain more alternations than the sequence $s_{i-1}, \ldots, s_k$, say $A'_2 > A_2$. But the number of alternations of $s_{n+1}, \ldots, s_k$ is at most equal to $A_1 + A_2 + 1$, because this sequence, being the concatenation of $s_n, \ldots, s_i$ and $s_{i-1}, \ldots, s_k$, contains only the alternations of these sequences, and possibly an additional alternation $s_i s_{i-1}$ (if $s_i \neq s_{i-1}$). Since $A'_1 + A'_2 \geq A_1 + 1 + A_2 + 1 > A_1 + A_2 + 1$, it follows that the sequence $s'_{n+1}, \ldots, s'_k$ contains more alternations than the sequence $s_n, \ldots, s_k$, in accordance to the induction hypothesis. ∎

Remarks:

1. Segner's lemma can be extended in the following manner:

*For any polynomial $Q$, multiplying $Q$ by $X - \alpha$, with $\alpha > 0$, increases the number of alternations of signs of $Q$ by an odd number, while multiplying it by $X + \alpha$ increases its number of permanences of signs by an odd number.*

This extension is obtained by including straightforward considerations of parity inside the induction above. But since it is a particular case of Fourier's rule below, it is essentially useless.

2. It should be clear that Segner's lemma and the above proof, and hence Descartes' rule, hold if the coefficients of the polynomials belong to any ordered field.

## 2.2. Fourier's Rule

We will now show that Fourier's rule is essentially independent of Descartes' rule, since each can be proved algebraically without the help of the other. The rule is:

*For any real polynomial $P$, $v(P)$ and $z^+(P)$ are of the same parity, as are $c(P)$ and $z^-(P)$.*

*Proof*: The second assertion follows immediately from the first by changing $P(X)$ into $P(-X)$.

We can assume without loss of generality that $P(0) \neq 0$, because if $P(0) = 0$, the polynomial $Q$ obtained by dividing $P$ by a suitable power $X^k$ fulfills $Q(0) \neq 0$ and contains the same number of alternations of signs and the same number of positive roots as $P$.

Let us write $P = A P_1 P_2$, where $P_1$ and $P_2$ are monic, $A \in \mathbb{R}$, $P_1$ has no real positive roots, and all the roots of $P_2$ are real positive.

Put $\deg(P_2) = n$. Since the coefficients of $P_1$ are real, its complex roots come in conjugate pairs $\gamma, \bar{\gamma}$. Therefore the polynomial $P_1$ can be written

$$P_1(X) = (X - \gamma_1)(X - \bar{\gamma}_1)(X - \gamma_2)(X - \bar{\gamma}_2) \cdots (X + \alpha_1)(X + \alpha_2) \cdots$$

where the $\gamma_i, \bar{\gamma}_i$ are the complex roots of $P_1$, and $-\alpha_i$ are its negative roots ($\alpha_i > 0$). The last coefficient of $P_1$ is clearly equal to

$$P_1(0) = (-\gamma_1)(-\bar{\gamma}_1)(-\gamma_2)(-\bar{\gamma}_2) \cdots \alpha_1 \alpha_2 \cdots = |\gamma_1|^2 |\gamma_2|^2 \cdots \alpha_1 \alpha_2 \cdots > 0.$$

Thus, the last coefficient of $P_1$ is strictly positive: $P_1(0) > 0$. Put

$$P_2 = (X - \beta_1)(X - \beta_2) \cdots,$$

where $\beta_1, \beta_2, \cdots > 0$. The last coefficient of $P$ is $P(0) = AP_1(0)P_2(0)$, hence is of the same sign as $AP_2(0) = A(-\beta_1)(-\beta_2)\cdots$. Consequently, the signs of $A$ and $P(0)$ are equal if and only if $P$ contains an even number of positive roots.

So, we see that the theorem is equivalent to the fact that a polynomial contains an even number of alternations of signs if and only if its leading and last coefficients have like signs. To this end, we need only show the following assertion:

CLAIM: *A sequence of signs contains an even number of alternations if and only if its extremities are equal.*

This is obviously true if the length of the sequence is 2. Assume inductively this assertion is true for sequences of lengths $n-1$, and consider a sequence of signs $S = (s_1, s_2, \ldots, s_n)$ of length $n$ ($s_i \in \{+, -\}$). Let $S' = (s_1, s_2, \ldots, s_{n-1})$. If $s_{n-1} = s_n$, then the number of alternations of $S$ and $S'$ are equal, and the extremities of $S$ and $S'$ are pairwise equal; since $S'$ fulfills the induction hypothesis, it is clear that $S$ fulfills it as well. If $s_{n-1} \neq s_n$, then the number of alternations of $S'$ is larger than the number of alternations of $S$ by 1; hence the parities of these numbers are opposite. But the extremities of $S$ are either opposite or equal, depending on whether the extremities of $S'$ are equal or opposite, resp. This implies that the induction hypothesis holds for $S$ in this case, as well. ∎

REMARK: The previous theorem and its proof hold obviously if the coefficients of $P$ belong to a real field. More generally, the above proof contains an interesting proposition, valid in any ordered field:

For every two polynomials $P$ and $Q$, $v(PQ) \equiv v(P) + v(Q) \pmod{2}$.

*Proof*: Let us denote by $A$ and $a$ the leading and last coefficient of $P$ resp., and by $B$ and $b$ the leading and last coefficient of $Q$ resp. Then the leading coefficient of $PQ$ is $AB$ and the last coefficient is $ab$. According to the claim in the above proof, $v(P) \equiv 0 \pmod{2}$ if and only if $A$ and $a$ share the same sign. Similarly, $v(Q) \equiv 0 \pmod{2}$ if and only if $B$ and $b$ share the same sign, and $v(PQ) \equiv 0 \pmod{2}$ if and only if $AB$ and $ab$ share the same sign. But since $(AB)(ab) = (Aa)(Bb)$, it is clear that the signs of $AB$ and $ab$ are equal if and only if $Aa$ and $Bb$ are both positive or both negative, that is, if $v(P)$ and $v(Q)$ are both even or both odd. Thus, $v(PQ)$ is even if and only if $v(P) + v(Q)$ is even. ∎

## 2.3. De Gua's Rule

De Gua's rule (*théorème des lacunes* in french textbooks) is:

*If, in a polynomial $P$, a group of $r$ consecutive terms is missing, then $P$ has at least $r$ imaginary roots if $r$ is even, or it has at least $r+1$ or $r-1$ imaginary roots if $r$ is odd, depending on whether the terms immediately preceding and following the group have like or unlike signs resp.*

It turns out that De Gua's rule is weaker than the following rule, which is an immediate corollary of Descartes' rule of signs:

**2.3.1. Corollary**: *The number of imaginary roots of a polynomial $P$ is at least equal to $\deg(P) - z^0(P) - v(P) - c(P)$.*

Therefore, all we must prove is that this corollary implies De Gua's rule.

*Proof*: In the process of counting the number of alternations and permanences of signs of $P$, only pairs of contiguous non-zero coefficients $(a_i, a_{i-1})$ are involved, or blocks of contiguous coefficients of the form $(a_i, 0, \ldots, 0, a_j)$ with $a_i a_j \neq 0$, containing one or more "0", say $r$ times "0". Now, a pair $(a_i, a_{i-1})$ contributes either one alternation, or one permanence of signs in the overall sum $v(P) + c(P)$. However, according to what has been explained at the beginning of Sec. 2, a block such as the one above contributes one alternation if $a_i a_j < 0$, 0 alternation if $a_i a_j > 0$, one permanence if $a_i a_j > 0$ and $r$ is even or if $a_i a_j < 0$ and $r$ is odd, and 0 permanence otherwise. Let $q$ be the number of consecutive pairs in the block: $q = r + 1$. If $r$ is even, we have just seen that the block contributes either one alternation, or one permanence. Therefore, the difference between the number of consecutive pairs and the sum of the numbers of alternations and permanences occurring in this block is $q - 1 = r$, a positive loss. Similarly, if $r$ is odd and $a_i a_j < 0$, then the block contributes one alternation and one permanence, hence the loss is $q - 2 = r - 1$. Finally, if $r$ is odd and $a_i a_j > 0$, then it contributes no alternation and no permanence, hence the loss is $q = r + 1$. In any case, the loss is $\geq 0$ (in the case of a normal contiguous pair $(a_i, a_{i-1})$ it is null). Furthermore, it is clear that the process of computing the loss associated with a block coincides exactly with De Gua's method. Since $\deg(P) - z^0(P)$ is equal to the number of contiguous pairs of non-trailing coefficients, it follows that $\deg(P) - z^0(P) - v(P) - c(P)$ is the sum of the losses of each block as above. This particularly implies De Gua's rule, and is even more powerful, as one can apply this rule to each sequence of contiguous zero (non-trailing) coefficients, summing the losses over the sequences. ∎

## 2.4. Budan's Rule and Fourier's theorem

As stated earlier in the paper, Budan's rule can be given a relatively simple algebraic proof based on Descartes' rule ([6]). Even so, it cannot be extended classes of functions more general than polynomials. Admittedly, the proof we present is not *ultra-simple*, but it has the merit of being valid for every function $f$ whose $n$-th differential does not vanish and is of constant sign inside an interval. In our opinion, a simpler proof of this result cannot be obtained without sacrificing rigor. We did not find this extension explicitly in the work of Fourier, but in regard to his proof, it is likely that he had an extension of this type in mind[10]. Actually, we believe the proof below is the exact expression of the ideas of Fourier. To state and prove this theorem, we need a few definitions and notations.

For every function f, we shall say that a function is differentiable inside a closed interval $I = [a, b]$ if it is differentiable in the interior of $I$, and is right differentiable at $a$ and left differentiable at $b$. From L'Hospital's rule, it follows that this condition is equivalent to the fact that $f$ is differentiable in the interior of $I$, and $f'$ has a finite limit at $a$ and $b$ (hence it can be extended by continuity at $a$ and $b$). In this form, this definition remains valid for intervals $[a, b] \in \overline{\mathbb{R}}$, with $a = -\infty$ and/or $b = +\infty$.

Given a point $\xi \in [a, b]$, we say that a property holds for every $x$ *sufficiently close* to $\xi$ (or *in the vicinity* of $x$), if there exists $\varepsilon > 0$ such that it holds for every $x \in [a, b]$ with $|x - \xi| < \varepsilon$. One could choose another metric, like

$$d(u, v) = |\arctan(u) - \arctan(v)|,$$

and this definition would be equivalent to the fact that $d(\xi, x) < \varepsilon$. Moreover, extending $\arctan(x)$ to $\overline{\mathbb{R}}$ in the obvious way, $d$ would extend, in turn, to a metric $\bar{d}$ over the compactified $\overline{\mathbb{R}}$, generating the usual topology of $\overline{\mathbb{R}}$. In this setting, the previous definition extends to every $\xi \in [a, b] \subseteq \overline{\mathbb{R}}$, replacing $|\xi - x| < \varepsilon$ by $\bar{d}(\xi, x) < \varepsilon$.

If $f$ is differentiable $k$ times inside $I$, then its $k$-th differential is denoted by $f^{(k)}$. By convention, $f$ is 0-differentiable if $f$ is continuous, and then $f^{(0)} = f$.

Given a function $f$ differentiable $n$ times, we denote by $V(f, t, n)$ the number of alternations of signs of the sequence

$$f^{(n)}(t), f^{(n-1)}(t), \ldots, f'(t), f(t) (= f^{(0)}(t)),$$

where we adopt the same minimization convention as we did above for sequences of polynomial coefficients (this amounts to ignore zero coefficients, and to count 0 alternations if the list contains no non zero, or a single non zero term). In the case where all the $f^{(i)}(t) = 0$, or if $n = 0$, we set $V(f, t, n) = 0$. It is obvious that $V(f, t, n) \in \mathbb{N}$, and $V(f, t, n) \leq n$.

If $f$ is a polynomial and $n = \deg(f)$, then $V(f, 0, n)$ is clearly the number of sign alternations of $f$ by Taylor's theorem. But a further property occurs: we have

$$V(f, +\infty, n) = 0 \quad \text{and} \quad V(f, -\infty, n) = n;$$

---

[10] Fourier explicitly stated that his theorem and the proof could be extended to more general classes of functions.

indeed, assuming for example $f^{(n)} > 0$, $f^{(n-1)}$ is increasing, hence $f^{(n-1)}(+\infty) = +\infty$. This implies that $f^{(n-1)} > 0$ near $+\infty$, hence $f^{(n-2)}(+\infty) = +\infty$ as well, and by an immediate induction, $f(+\infty) = f'(+\infty) = \ldots = f^{(n)}(+\infty) > 0$. A similar argument shows that if $f^{(n)} < 0$, $f(+\infty) = f'(+\infty) = \ldots = f^{(n)}(+\infty) < 0$. In both cases, the sequence of differentials has no alternations of signs, that is, $V(f, +\infty, n) = 0$.

If now we set $f(x) = g(-x)$, we have $f^{(i)}(x) = (-1)^i g^{(n)}(-x)$, hence, by what we have just proved, $\text{sign}(f^{(i)}(-\infty)) = (-1)^i$. This shows that $V(f, -\infty, n) = n$. ∎

We say that a zero $\alpha$ of $f$ is of multiplicity $\mu$ if $f^{(i)}(\alpha) = 0$ for all $0 \leq i \leq \mu$, and $f^{\mu+1}(\alpha) \neq 0$. By convention, $\alpha$ is of multiplicity 0 if $\alpha$ is not a zero of $f$. We also use $Z(f, I)$ to denote the number of zeros of $f$ inside an interval $I$, counted with their multiplicities (by convention, $Z(f, ]\alpha, \alpha]) = Z(f, \emptyset) = 0$).

Finally, given any number $x \in \overline{\mathbb{R}}$, we note $\text{sign}(x) = 1$ if $x$ is positive, $\text{sign}(x) = -1$ if $x$ is negative, and $\text{sign}(x) = 0$ if $x = 0$.

Budan's rule is:

If $P$ is a non-zero polynomial of degree $n$, and $]a, b] \subseteq \mathbb{R}$, then the number of roots of $P$ inside $]a, b]$ is at most equal to

$$V(P, a, n) - V(P, b, n).$$

Furthermore, equality holds for every $a < b$ if, and only if, all the roots of $P$ are real (Loewy, Curtis).

The second assertion follows easily from the first, by observing that if $a \leq c \leq b$,

$$V(P, a, n) - V(P, b, n) = V(P, a, n) - V(P, c, n) + V(P, c, n) - V(P, b, n),$$
$$Z(P, ]a, b]) = Z(P, ]a, c]) + Z(P, ]c, b]),$$
$$V(P, -\infty, n) = n,$$
$$V(P, +\infty, n) = 0.$$

Regarding the first assertion, since the $n$th differential of $P$ is constant and $\neq 0$, it is a particular case of the following theorem, essentially due to Fourier:

**2.4.1. Theorem** (Fourier's theorem): *Let $[a, b] \subseteq \mathbb{R}$, and $f \colon [a, b] \to \mathbb{R}$ be $n$ times differentiable inside $[a, b]$. Assume that the $n$-th derivative $f^{(n)}$ does not vanish and is of constant sign inside $[a, b]$. Then*

*(i) the function $x \mapsto V(f, x, m)$ is decreasing and right continuous inside $[a, b]$, for every $m \leq n$. In particular, $V(f, a, n) - V(f, b, n) \geq 0$.*

*(ii) the number of zeros of $f$ inside $]a, b]$, including multiplicities, is not larger than*

$$V(f, a, n) - V(f, b, n),$$

*and these two numbers share the same parity.*

*Proof*: For the sake of simplicity, we shall say that $\mathcal{T}(u,v,n)$ is true if $a \leq u < v \leq b$ and if

$$Z(f,]u,v]) = V(f,u,n) - V(f,v,n) - 2s, \quad \text{with} \quad s \in \mathbb{N}.$$

If $n$ has been fixed, we also abbreviate $\mathcal{T}(u,v,n)$ by $\mathcal{T}(u,v)$. It should be clear that $\mathcal{T}$ is additive in the following sense: If $\mathcal{T}(u,v,n)$ and $\mathcal{T}(v,w,n)$ are true, then $\mathcal{T}(u,w,n)$ is true. This can be seen with ease, adding term by term the equations corresponding to $\mathcal{T}(u,v,n)$ and $\mathcal{T}(v,w,n)$.

Now, it is clear that the second assertion of the theorem is equivalent to the fact that $\mathcal{T}(a,b,n)$ is true for all $n \in \mathbb{N}$. If $n = 0$, then $Z(f,]a,b]) = 0$ and $V(f,x,n) = 0$ for all $x$, hence the correctness of the theorem is obvious in this case.

The following lemma shows that we must only prove the theorem locally. It is a bit more general that what we need presently, but this generality will be useful later.

**2.4.2. Lemma**: *Let $[a,b] \subseteq \overline{\mathbb{R}}$, and $\mathcal{T}(u,v)$ be a property depending on $u$ and $v$, with $a \leq u \leq v \leq b$. Assume that $\mathcal{T}$ is additive in the sense that for all points $\xi$ of $[a,b]$, the condition*

"*$u \leq \xi \leq v$ and both $\mathcal{T}(u,\xi)$ and $\mathcal{T}(\xi,v)$ hold*"

*implies that $\mathcal{T}(u,v)$ holds.*

*Then the following conditions are sufficient in order to ensure that $\mathcal{T}(a,b)$ holds:*

- *For every $\xi \in ]a,b[$, except, perhaps, a given set $S$ of isolated points not containing $a$ and $b$, $\mathcal{T}(u,v)$ hold whenever $u$ and $v$ are sufficiently close to $\xi$ and fulfill $a \leq u \leq \xi \leq v \leq b$;*

- *at those points $\xi \in S$ (if $S$ is not empty), $\mathcal{T}(u,v)$ holds for every $u,v$ sufficiently close to $\xi$ and such that $a \leq u < \xi < v \leq b$.*

*Proof*: According to the hypothesis, $\mathcal{T}(a,v)$ holds for every $v$ sufficiently close to $a$. Let $\xi$ be the upper bound of all the numbers $x$ such that $\mathcal{T}(a,t)$ holds for every $t \leq x$, at the possible exception of those $t$ belonging to the set $S$ of isolated points specified in the lemma.

Assume, in order to obtain a contradiction, that $\xi \neq b$. It is clear that $\mathcal{T}(a,x)$ holds for every $x < \xi$ with $x \notin S$.

According to the hypothesis, for every $u, v$ sufficiently close to $\xi$, with $a \leq u < \xi < v \leq b$, $\mathcal{T}(u,v)$ holds. Reducing eventually the vicinity zone of $\xi$, it can be supposed that $[u,v]$ contains no element of $S$, except, perhaps, $\xi$. Because of the definition of $\xi$, $\mathcal{T}(a,u)$ holds, and the same is true of $\mathcal{T}(u,x)$ for every $x$ with $\xi < x \leq v$. By the additivity of $\mathcal{T}$, it follows that $\mathcal{T}(a,x)$ holds for every $x$ with $\xi < x \leq v$. If now $x = \xi$ and $\xi \notin S$, then the strict inequalities above can be replaced by non strict ones (hypothesis), from what follows that $\mathcal{T}(a,x)$ holds as well. This shows that $\mathcal{T}(a,x)$ holds for every $x \notin S$ such that $x \leq v$, in contradiction with the fact that $\xi$ is the upper bound of all such numbers $v$. Thus $\xi = b$.

This does not prove yet that $\mathcal{T}(a,b)$ holds, but since $b \notin S$ (hypothesis), one can use the additivity of $\mathcal{T}$ and the fact that $\mathcal{T}([u,b])$ holds for every $u$ in the vicinity of $b = \xi$, to show, exactly as above, that $\mathcal{T}(a,b)$ holds. ∎

Considering the above lemma, it suffices to prove that if $\xi \in [a,b]$, the properties $\mathcal{T}(u,v) := \mathcal{T}(u,v,n)$ are true for every $u,v$ sufficiently close to $\xi$, and $a \leq u \leq \xi \leq v \leq b$. At the same time, we prove that $V(f,x,m)$ is right continuous at $\xi$ as a function of the variable $x$, for every $\xi \in [a,b]$ and $m \leq n$.

These assertions being correct if $n = 0$, let us assume inductively that they are for every $n \leq N$, with $N \in \mathbb{N}$. Let $n = N+1$, $0 < m \leq n$, and suppose, as in the theorem, that $f^{(n)}$ keeps its signs constant inside $[a,b]$.

Assume first that for some $k$ with $0 < k < n$, $f^{(k)}(\xi) \neq 0$. In this case, $f^{(k)}(x)$, being differentiable, is continuous at $\xi$, and hence does not vanish and does not change its sign in the vicinity of $\xi$. Consequently, the induction hypothesis (with $f^{(k)}$ in place of $f^{(n)}$ or $f^{(k)}$ in place of $f$) implies that for every $u,v$ in the vicinity of $\xi$, with $a \leq u \leq \xi \leq v \leq b$,

$$Z(f,]u,v]) = V(f,u,k) - V(f,v,k) - 2s$$
$$\text{and} \quad Z(f^{(k)},]u,v]) = V(f^{(k)},u,n-k) - V(f^{(k)},v,n-k) - 2s'.$$

But $Z(f^{(k)},]u,v]) = 0$ since $f^{(k)}$ does not vanish in $[u,v]$. Moreover, it is clear in general, that if $f^{(k)}(\alpha) \neq 0$,

$$V(f,\alpha,n) = V(f,\alpha,k) + V(f^{(k)},\alpha,n-k).$$

Therefore, adding the two previous equalities term by term, there holds

$$Z(f,]u,v]) = V(f,u,n) - V(f,v,n) - 2s - 2s' = V(f,u,n) - V(f,v,n) - 2s''.$$

Thus, $\mathcal{T}(u,v)$ is true in the case where $f^{(k)}(\xi) \neq 0$ for some $0 < k < n$, for every $u,v$ in the vicinity of $\xi$, with $a \leq u \leq \xi \leq v \leq b$.

As a particular case, if $f^{(k)}(a) \neq 0$, then $\mathcal{T}(a,v)$ holds for all $v$ sufficiently close to $a$.

Furthermore the induction hypothesis implies that $V(f,x,k)$ and $V(f^{(k)},x,n-k)$ are right continuous at $\xi$, as functions of $x$, hence, so is

$$V(f,\xi,n) = V(f,\xi,k) + V(f^{(k)},\xi,n-k).$$

The induction step has been performed in the case where $f^{(k)}(\xi) \neq 0$ for some $k \in \{1, \ldots, n-1\}$. Henceforth, we suppose that $f^{(k)}(\xi) = 0$, for every $k$ with $0 < k < n$. In order to show that $\mathcal{T}(u,v)$ holds for every $u, v$ in the vicinity of $\xi$ ($a \leq u \leq \xi \leq v \leq b$), it suffices to show, because of the additivity of $\mathcal{T}$, that $\mathcal{T}(u, \xi)$ and $\mathcal{T}(\xi, v)$ hold.

We first prove that $\mathcal{T}(\xi, v)$ holds for all $v > \xi$ sufficiently close to $\xi$, and consider, to this end, the signs of $f^{(k)}(x)$ in the vicinity and to the right of $\xi$. By the mean value theorem, $f^{(n-1)}$ is strictly increasing or strictly decreasing inside $[a, b]$, according to whether $f^{(n)}$ is positive or negative inside $]a, b[$. But $f^{(n-1)}(\xi) = 0$, hence

$$\text{sign}(f^{(n-1)}) = \text{sign}(f^{(n)}) \text{ inside } ]\xi, b].$$

It follows that $f^{(n-2)}$ is strictly increasing or decreasing inside $[\xi, b]$, according to whether $f^{(n-1)}$ is positive or negative. Since $f^{(n-2)}(\xi) = 0$,

$$\text{sign}(f^{(n-2)}) = \text{sign}(f^{(n-1)}) = \text{sign}(f^{(n)}) \text{ inside } ]\xi, b].$$

This induction can be pursued until $f^{(1)} = f'$ is reached (or even $f^{(0)} = f$ if $f(\xi) = 0$), showing that

$$\text{sign}(f^{(k)}(x)) = \text{sign}(f^{(n)}(x)), \text{ for all } 0 < k < n \text{ and } x \in ]\xi, b],$$

and this relation continues to hold anyway for $k = 0$, whenever $f(\xi) = 0$. So, $V(f, v, n)$, or the number of alternations of signs inside the sequence

$$f(v), f'(v'), f''(v), \ldots, f^{(n-1)}(v), f^{(n)}(v),$$

is the same as the number of alternations of signs in the sequence $f(v), f^{(n)}(v)$. On the other hand, it is immediate that $V(f, \xi, n)$ is equal to the number of alternations of signs in the sequence $f(\xi), f^{(n)}(\xi)$ (since $f'(\xi) = \ldots = f^{(n-1)}(\xi) = 0$). Moreover, we have:

1) $\text{sign}(f^{(n)}(v)) = \text{sign}(f^{(n)}(\xi)$ by hypothesis.

2) If $f(\xi) \neq 0$, then $\text{sign}(f(v)) = \text{sign}(f(\xi))$ for every $v \geq \xi$ sufficiently close to $\xi$ (continuity of $f$).

3) If $f(\xi) = 0$, then $\text{sign}(f(v)) = \text{sign}(f^{(n)}(v))$ for every $v > \xi$ sufficiently close to $\xi$, as shown just above.

From these considerations, it follows easily that

$$V(f, \xi, n) - V(f, v, n) = 0,$$

whether or not $f(\xi) = 0$. Hence, $V(f, x, n)$ is right continuous at $\xi$ as a function of the variable $x$, ending the proof of this assertion for every $\xi \in [a, b]$.

In addition, since $f(v) \neq 0$ for every $v > \xi$ sufficiently close to $\xi$, it is clear that $Z(f, ]\xi, v]) = 0$. Therefore $\mathcal{T}(\xi, v)$ is true for every $v > \xi$ sufficiently close to $\xi$, as was to be shown. In particular, $\mathcal{T}(a, v)$ holds for every $v > a$ in the vicinity of $a$.

Now, in order to show that $\mathcal{T}(u, \xi)$ holds for all $u < \xi$ sufficiently close to $\xi$, we consider the sign of $f^{(k)}(x)$ in the vicinity of $\xi$, but this time on the left side of $\xi$. By the mean value theorem, $f^{(n-1)}$ is strictly increasing or strictly decreasing inside $[a, b]$, according to whether $f^{(n)}$ is positive or negative inside $]a, b[$. But $f^{(n-1)}(\xi) = 0$, hence

$$\operatorname{sign}(f^{(n-1)}) = -\operatorname{sign}(f^{(n)}) \text{ inside } [a, \xi[ \ .$$

It follows that $f^{(n-2)}$ is strictly increasing or decreasing inside $[a, \xi[$, according to whether $f^{(n-1)}$ is positive or negative. Since $f^{(n-2)}(\xi) = 0$,

$$\operatorname{sign}(f^{(n-2)}) = -\operatorname{sign}(f^{(n-1)}) = \operatorname{sign}(f^{(n)}) \text{ inside } [a, \xi[ \ .$$

This induction can be pursued until $f^{(1)} = f'$ is reached (or even $f^{(0)} = f$ if $f(\xi) = 0$), showing that

$$\operatorname{sign}(f^{(k)}(x)) = (-1)^{n-k} \operatorname{sign}(f^{(n)}(x)), \text{ for every } 0 < k < n \text{ and } x \in [a, \xi[ \ .$$

Moreover, this relation holds for $k = 0$, whenever $f(\xi) = 0$. Hence, if $f(\xi) = 0$ and $u$ is sufficiently close to $\xi$, with $a \leq u < \xi$, then

$$V(f, u, n) = n \text{ and } Z(f, ]u, \xi]) = n.$$

On the other hand, it is evident that $V(f, \xi, n) = 0$, since $f^{(k)}(\xi) = 0$ for every $k < n$, therefore $\mathcal{T}(u, \xi)$ holds if $f(\xi) = 0$.

If, in contrast, $f(\xi) \neq 0$ and $a \leq u < \xi$, then whenever $u$ is sufficiently close to $\xi$, the continuity of $f$ implies that $f$ does not vanish inside $[u, \xi]$, therefore

$$Z(f, ]u, \xi]) = 0 \text{ and } \operatorname{sign}(f(u)) = \operatorname{sign}(f(\xi)).$$

But $V(f, u, n)$, which is the number of alternations of signs inside the sequence

$$f(u), f'(u), \ldots, f^{(n)}(u),$$

is even if and only if the extremities $f(u)$ and $f^{(n)}(u)$ share the same sign: this can be shown by an easy induction (see Sec. 2.2), or more intuitively, considering the fact that someone who has to cross a river a certain number of times, must cross it an even number of times if (and only if) he wishes to return to the side of the river from which he has started from. On the other hand, it is clear that $V(f, \xi, n)$ is equal to 0 if $f(\xi)$ and $f^{(n)}(\xi)$ share the same sign, and to 1 otherwise (since $f^{(k)}(\xi) = 0$ for all $0 < k < n$). Thus, since the signs of $f(\xi)$ and of $f(u)$ are identical, as well as the signs of $f^{(n)}(\xi)$ and $f^{(n)}(u)$ (hypothesis), we conclude that $V(f, u, n) - V(f, \xi, n)$ is always a non-negative even number, and we have

$$Z(f, ]u, \xi]) = 0 = V(f, u, n) - V(f, \xi, n) - 2s, \quad \text{with} \quad s \in \mathbb{N}.$$

This ends the proof that $\mathcal{T}(u, v)$ holds in any case, whenever $u$ and $v$ are sufficiently close to $\xi$, and $a \leq u \leq \xi \leq v \leq b$.

The demonstration that $\mathcal{T}(a, b)$ holds is now achieved.

Finally, exchanging $a$ and $b$ with $x$ and $x'$ resp. in the result above, we have

$$\begin{aligned} Z(f, ]x, x']) &= V(f, x, n) - V(f, x', n) - 2s \\ \text{or} \quad V(f, x, n) - V(f, x', n) &= Z(f, ]x, x']) + 2s \geq 0. \end{aligned}$$

for every $x$ and $x' > x$. Hence $V(f, x, n)$ is decreasing, showing the first assertion of the theorem. ∎

REMARK: Given a function $f$ $m$ times differentiable inside $[a, b]$, it may be asked whether $V(f, a, m) - V(f, b, m)$ can be negative. Actually, it can, as shown by the function
$$f(x) = x^2 - x - 1, \quad \text{with} \quad [a, b] = [0, 1], \quad \text{and} \quad m = 1.$$
Indeed, $f(0) = -1$, $f'(0) = -1$, $f(1) = -1$ and $f'(1) = 1$.

## 2.5. Hurwitz's theorem

Fourier's theorem (second assertion) is, in turn, a particular case of the following proposition, that is a variation of a theorem due to Hurwitz ([22]). It will be given its final form in Theorem 2.5.2

**2.5.1. Proposition**: *Using the definitions, notations and conventions stated in the previous section, let us assume that $[a, b] \subseteq \overline{\mathbb{R}}$ and that $f \colon [a, b] \to \mathbb{R}$ is $n$ times differentiable inside $[a, b]$ ($n \in \mathbb{N}$). Furthermore, let us assume that the $n$-th derivative $f^{(n)}$ does not vanish and is of constant sign inside $[a, b]$.*

*If, for some $m \in \mathbb{N}$ with $0 \leq m \leq n$, $f^{(m)}(a) \neq 0$ and $f^{(m)}(b) \neq 0$, then*

$$Z(f, ]a, b]) = Z(f^{(m)}, ]a, b]) + V(f, a, m) - V(f, b, m) - 2s, \quad \text{with } s \in \mathbb{N}.$$

REMARKS:

1. If $m = n$, this theorem is nothing else than Fourier's theorem.

2. As pointed out at the end of the previous section, nothing prevents $V(f,a,m) - V(f,b,m)$ from being negative.

3. The assumption that $f^{(m)}(b) \neq 0$ is essential, as shown by the following example: Let $f(x) = (x-1)^3 - 1$. Then $Z(f, ]0,1]) = 0$, and $Z(f'', ]0,1]) = 1$. moreover, the sequence
$$f(0) \quad f'(0) \quad f''(0) \quad \text{is} \quad -2 \quad 3 \quad -6,$$
and the sequence
$$f(1) \quad f'(1) \quad f''(1) \quad \text{is} \quad -1 \quad 0 \quad 0.$$
Hence $V(f,0,2) - V(f,1,2) = 2$, and it is impossible that
$$Z(f, ]0,1]) = 0 = Z(f'', ]0,1]) + V(f,0,2) - V(f,1,2) - 2s = 3 - 2s.$$

Similarly, if $f(x) = -x^2 + 2$, $m = 1$ and $[a,b] = [0,1]$, Then $V(f,0,1) = 0$ and $V(f,1,1) = 1$. On the other hand, $Z(f, ]a,b]) = Z(f', ]a,b]) = 0$. This shows, as above, that the assumption $f^{(m)}(a) \neq 0$ is essential.

*Proof of the theorem*: Reusing some parts of the proof of Fourier's theorem above, it would have been possible to give a direct proof of this result. Rather, in order not to repeat ourselves, we shall use Fourier's theorem inside the proof.

Given a number $m$ with $0 \leq m \leq n$, we say that $\mathcal{T}(u,v,m)$ is true if $a \leq u < v \leq b$ and if
$$Z(f, ]u,v]) = Z(f^{(m)}, ]u,v]) + V(f,u,m) - V(f,v,m) - 2s, \quad \text{with} \quad s \in \mathbb{N}.$$

If $m$ has been fixed, we also abbreviate $\mathcal{T}(u,v,m)$ by $\mathcal{T}(u,v)$. As previously, it should be clear that $\mathcal{T}$ is additive in the following sense: If $\mathcal{T}(u,v,m)$ and $\mathcal{T}(v,w,m)$ are true, then $\mathcal{T}(u,w,n)$ is true.

We have to prove that $\mathcal{T}(a,b,m)$ is true for all $n \in \mathbb{N}$ whenever $f^{(n)}$ does not vanish, is of constant sign in $[a,b]$, and $f^{(m)}(a)$ and $f^{(m)}(b)$ are not equal to 0 for some $m \leq n$. If $m = n$, this assertion coincides with the second assertion of Fourier's theorem. Hence, the theorem is true in this case. In particular, it is true if $m = n = 0$.

Given $\xi \in [a,b]$ and $m \leq n$, we first prove that $\mathcal{T}(u,v) = \mathcal{T}(u,v,m)$ is true, for every $u, v$ sufficiently close to $\xi$, with $f^{(m)}(u) f^{(m)}(v) \neq 0$. Since this assertion is true for $m = n = 0$, let us assume inductively that it is true for every $m \leq n \leq N$, with $N \in \mathbb{N}$. Let $n = N + 1$, and suppose, as in the theorem, that $f^{(n)}$ does not change its sign inside $]a,b]$. As noted above, the theorem is true if $m = n$, hence we can assume from now on that $m < n$.

If, for some $k$ with $m \leq k < n$, $f^{(k)}(\xi) \neq 0$, then the sign of $f^{(k)}$ must remain constant in the vicinity of $\xi$, because of the continuity of $f^{(k)}$. Hence, for every $u$ and $v$ sufficiently close to $\xi$, the induction hypothesis (with $k$ in place of $n$) immediately implies that $\mathcal{T}(u,v,m)$ is true. In particular, $\mathcal{T}(u,v,m)$ holds for $\xi = b = v$, because of the hypothesis in the theorem.

Henceforth, we now suppose that $f^{(k)}(\xi) = 0$ for every $k \geq m$, and that $\xi < b$.

If, for some $k$ with $0 < k < m$, holds $f^{(k)}(\xi) \neq 0$, then $f^{(k)}$ does not change its sign in the vicinity of $\xi$, as explained above. Therefore, for all $u,v$ sufficiently close to $\xi$, the induction hypothesis implies that

$$Z(f, ]u,v]) = Z(f^{(k)}, ]u,v]) + V(f,u,k) - V(f,v,k) - 2s$$
$$\text{and} \quad Z(f^{(k)}, ]u,v]) = Z(f^{(m)}(]u,v]) + V(f^{(k)}, u, m-k) - V(f^{(k)}, v, m-k) - 2s'.$$

But in general, if $f^{(k)}(\alpha) \neq 0$, it is plain that

$$V(f, \alpha, n) = V(f, \alpha, k) + V(f^{(k)}, \alpha, n-k), \tag{4}$$

as can be immediately seen. Hence, the two equations above imply

$$Z(f, ]u,v]) = Z(f^{(m)}, ]u,v]) + V(f,u,m) - V(f,v,m) - 2s'',$$

with $s'' = s + s'$. This shows that $\mathcal{T}(u,v)$ is true in this case.

It remains to consider the case where $f^{(k)}(\xi) = 0$, for every $k$ with $0 < k < n$, with $a \leq \xi < b$. Let us assume that this holds. Notice that since $f^{(m)}(a) \neq 0$ by hypothesis, $\xi \neq a$. Furthermore, Fourier's theorem implies that the number of zeros of $f^{(m)}$ is finite inside $]a,b]$ (since $V(f^{(m)}, a, n-m) - V(f^{(m)}, b, n-m)$ is at most equal to $n-m$). Let $S$ denotes the set of zeros of $f^{(m)}$ in $[a,b]$: $\xi \in S$. Since $S$ is finite, its points are isolated, hence $f^{(m)}(x) \neq 0$ for every $x \neq \xi$ sufficiently close to $\xi$. It follows that for every $u,v$ in the vicinity of $\xi$, with $a < u < \xi < v < b$, $f^{(m)}(u) \neq 0$ and $f^{(m)}(v) \neq 0$. Also, $Z(f^{(m)}, ]u,v]) = n - m$ since $f^{(k)}(\xi) = 0$ for every $m \leq k < n$ (that is, the unique zero $\xi$ of $f^{(m)}$ in $[u,v]$ is of multiplicity $n-m$). Fourier's theorem can be applied between $u$ and $v$:

$$Z(f, ]u,v]) = V(f,u,n) - V(f,v,n) - 2s,$$
$$\text{and} \quad Z(f^{(m)}, ]u,v]) = V(f^{(m)}, u, n-m) - V(f^{(m)}, v, n-m) - 2s'.$$

Actually, since $Z(f^{(m)}, ]u,v]) = n-m$, and since $V(f^{(m)}, u, n-m) - V(f^{(m)}, v, n-m)$ is at most equal to $n-m$, it is necessary that $s' = 0$ in order for the second equation above to hold. So, subtracting term by term the second equation from the first one and taking (4) into account, we have

$$Z(f, ]u,v]) = Z(f^{(m)}, ]u,v]) + V(f,u,m) - V(f,v,m) - 2s.$$

This ends the proof of the fact that $\mathcal{T}(u,v)$ is true whenever $f^{(m)}(u)f^{(m)}(v) \neq 0$, and one of the two conditions hold:

"$\xi \notin S$ and $a \leq u \leq \xi \leq v \leq b$, with $u,v$ in the vicinity of $\xi$", or "$\xi \in S$ and $a \leq u < \xi < v \leq b$, with $u,v$ in the vicinity of $\xi$".

But since the points of $S$ are isolated, the vicinity zone around $\xi$ can be eventually reduced in such a way that $f^{(m)}(u)f^{(m)}(v) \neq 0$ for all $u,v$ sufficiently close to $\xi$, with $u,v \neq \xi$ if $\xi \in S$. Therefore the condition $f^{(m)}(u)f^{(m)}(v) \neq 0$ can be dropped in the assertion above; by Lemma 2.4.2, we can now conclude that $\mathcal{T}(a,b)$ holds, as was to be shown. ∎

There exists a more general result, essentially due to Hurwitz ([22]) in the case of regular functions. The proof given here is different from the proof of Hurwitz, and may be more easy. Also, Hurwitz stated and proved this theorem for regular functions $f$ only, but it is susceptible of a much greater extent. To this end, we need an additional definition.

DEFINITION: Assume that $f$ is a real function, and that $k \in \mathbb{N} \cup \infty$. We say that $f$ is *Taylor analyzable of order $k$* inside an interval $I$, if $f$ is at least $k$ times differentiable inside $I$, and if one of the two following conditions holds:

P1. $k < \infty$, and for every $\xi \in I$, there exists $n \geq k$ such that $f$ is $n$-times differentiable in $[a,b]$ and $\text{sign}(f^{(n)}(x))$ is constant and not null, for every $x \in I$ in the vicinity of $\xi$.

P2. $k = \infty$, and for every $\xi \in I$ and every $m \in \mathbb{N}$, there exists $n \geq m$ such that $\text{sign}(f^{(n)}(x))$ is constant and not null, for every $x \in I$ in the vicinity of $\xi$.

In the later case, $f^{(n)}$ is continuous for every $n \in \mathbb{N}$ inside $I$ (since it is differentiable), hence an equivalent condition is

P2'. $k = \infty$ and for every $\xi \in I$ and every $m \in \mathbb{N}$, there exists $n \geq m$ such that $f^{(n)}(\xi) \neq 0$.

It is evident that if $f$ is Taylor analyzable of order $k$, then $f^{(i)}$ is Taylor analyzable of order $k-i$, for every $i \leq k$. Also, we point out that if $f$ is Taylor analyzable of infinite order inside $I$ and if $I$ is compact, then $f$ is in fact Taylor analyzable of finite order $k$, for every $k \in \mathbb{N}$.

Typically, $f$ is Taylor analyzable of order $k$ inside $I$ in the following cases:

1) $k < \infty$, $f$ is at least $k'$ times differentiable inside $I$, with $k' \geq k$, $f^{(k')}$ is continuous, and for every $x \in I$, there exists $k \leq n \leq k'$ such that $f^{(n)}(x) \neq 0$.

2) $f$ is analytic inside $I$.

Indeed, in the first case, the continuity of $f^{(n)}$ implies that the sign of $f^{(n)}$ does not change in the vicinity of an element $\xi$, which is such that $f^{(n)}(\xi) \neq 0$. In the second case, the theorem of analytic continuation implies that $f$ either reduces to a polynomial, in which case it is obviously Taylor analyzable of finite order, or satisfies condition $P2'$.

The following theorem can now be stated and proved.

**2.5.2. Theorem**: *Let $[a,b] \subseteq \overline{\mathbb{R}}$ and $f \colon [a,b] \to \mathbb{R}$ be Taylor analyzable of order $k$ inside $[a,b]$, with $k \in \mathbb{N} \cup \infty$. Then*
*(i) the number of zeros of $f$ inside $[a,b]$ is finite;*
*(ii) if, for some $m \in \mathbb{N}$ with $0 \leq m \leq k$, $f^{(m)}(a) \neq 0$ and $f^{(m)}(b) \neq 0$, then*
$Z(f, ]a,b]) = Z(f^{(m)}, ]a,b]) + V(f,a,m) - V(f,b,m) - 2s$, *with* $s \in \mathbb{N}$.

*Proof*: Again, a direct proof could be given, reusing some part of the proof of Prop. 2.5.1. Rather, we shall use Prop. 2.5.1 inside the present proof.

(i) Since $f$ is Taylor analyzable of order $k$, for every $\xi \in [a,b]$, there exists $n$ such that $f^{(n)}$ does not vanish and is of constant sign in the vicinity of $\xi$. By Taylor's theorem, for every $\xi \in [a,b]$ and $u \leq \xi \leq v$, the number of zeros of $f$ inside $I_\xi = [u,v]$ is at most equal to $n$ whenever $u$ and $v$ are sufficiently close to $\xi$. Since $[a,b]$ is compact, it can be covered by a finite number of intervals $I_\xi$ of this form, say $I_{\xi_1}, I_{\xi_2}, \ldots$. Of course, the number of zeros of $f$ inside $[a,b]$ is at most equal to the sum of the number of zeros of $f$ inside the $I_{\xi_i}$. Thus, $Z(f,[a,b])$ is finite.

(ii) As previously, we say that $\mathcal{T}(u,v,m)$ is true if

$$Z(f, ]u,v]) = Z(f^{(m)}, ]u,v]) + V(f,u,m) - V(f,v,m) - 2s, \text{ with } s \in \mathbb{N}.$$

And again, $\mathcal{T}$ is additive: if $\mathcal{T}(u,v,m)$ and $\mathcal{T}(v,w,m)$ hold, then $\mathcal{T}(u,w,m)$ holds. We have to prove that $\mathcal{T}(a,b,m)$ holds. Assume that $\xi \in [a,b]$. It is immediate that $f^{(m)}$ is Taylor analyzable of order $k-m$. According to (i), it follows that $f^{(m)}$ has finitely many zeros inside $[a,b]$, and these zeros are isolated. Let $S$ denotes the set of zeros of $f^{(m)}$ in $[a,b]$. Notice that $a,b \notin S$ (hypothesis). Since the points of $S$ are isolated, for every $\xi \in ]a,b[$ and $u,v$ sufficiently close to $\xi$, with $a \leq u < \xi < v \leq b$, $f^{(m)}(u) \neq 0$ and $f^{(m)}(v) \neq 0$. Moreover, the previous strict inequalities can be replaced by wide ones if $\xi \notin S$, and the condition $\xi \in ]a,b[$ can be replaced by $\xi \in [a,b]$ as well. On the other hand, since $f$ is Taylor analyzable of order $k$, there exists $n \geq m$ such that the sign of $f^{(n)}$ is constant in the vicinity of $\xi$. From Prop. 2.5.1, it follows that $\mathcal{T}(u,v,m)$ holds if one of the following two conditions hold:

"$\xi \notin S$, $a \leq u \leq \xi \leq v \leq b$ with $u,v$ in the vicinity of $\xi$", or "$\xi \in S$ and $a \leq u < \xi < v \leq b$ with $u,v$ in the vicinity of $\xi$".

Thus, to end the proof of the theorem, it suffices to apply Lemma 2.4.2. ∎

Acknowledgments: The author wishes to thank Jacques Gélinas for his helpful informations and suggestions.

# Bibliography

[1] A. A. Albert. An inductive proof of Descartes's rule of sign. *The American Mathematical Monthly*, 50(3):178–180, Mar. 1943.

[2] R. D. Arthan. Descartes Rule of Signs by an Easy Induction. *Arxiv*, Oct.

[3] Margherita Bartolozzi and Raffaella Franci. La regola dei segni dall'enunciato di R. Descartes (1637) alla dimostrazione di C. F. Gauss (1828). (Italian) [The rule of signs enunciated by R. Descartes (1637) as demonstrated by C. F. Gauss (1828)]. 45(4):335–374, December 1993.

[4] F. D. Budan. *Nouvelle méthode pour la résolution des equations d'un degré quelconque*. Bondé-Dupré, Bachelier, Paris, 1822.

[5] E. Rouché C. Hermite, H. Poincaré. *Œuvres de Laguerre*, volume 1. Gauthier-Villars et fils, Paris, 1898.

[6] N. B. Conkwright. An elementary proof of the Budan-Fourier Theorem. *The American Mathematical Monthly*, 50(10):603–605, Dec. 1943.

[7] D. R. Curtis. Extensions of Descartes' rule of signs connected with a problem suggested by Laguerre. *Trans. Amer. Math. Soc.*, 16:350–360, 1915.

[8] D. R. Curtis. Recent extensions of Descartes's Rule of Signs. *The Annals of Mathematics*, 19:251–278, Jun. 1918.

[9] G .Darboux. *Œuvres de Fourier*, volume 2, page 310. Gauthier-Villars & fils, Paris, 1890.

[10] R. Descartes. *La géometrie (Discours de la méthode, third part)*, page 373. Ed. of Leyde, 1637.

[11] R. Descartes. *Lettres de Mr Descartes*, volume 3. Charles Angot, Paris, 1667.

[12] L. E. Dickson. *First Course in the Theory of Equations*, pages 71–74 and 83–85. John Willey & sons, New York, 1922.

[13] Encyclopæ dia Britannica. Equations. In *Supplement to the Fourth, Fifth, and Sixth editions of the Encyclopædia Britannica*, volume 4, pages 674–675. Edinburgh, 1824.

[14] G. Eneström. In *Bibliotheca Mathematica*, page 307. Druck Und Verlag Von B. G Teubner, Leipzig, 1906. https://archive.org/stream/bibliothecamath01enesgoog#page/n14/mode/1up.

[15] J. Fourier. Sur l'usage du théorème de Descartes dans la recherches des limites des racines. *Bulletin des Sciences par la Société philomatique de Paris*, pages 156–165 and 181–187, Oct. and Dec. 1820.

[16] B. Gagneux. La règle des signes de Descartes : le long cheminement d'une imprécision. In *Mathématiciens francais du XVIIe siècle, Descartes, Fermat, Pascal*, pages 136–163. Presses Universitaires Blaise Pascal, April 2008.

[17] C. F. Gauss. Beweis eines algebraischen Lehrsatzes. *Crelle's Journal fur die reine und angewandte Mathematik*, 3, Jan. 1828.

[18] J. A. Grunert. Beweis des Harriotischen Satzes. *Journal für die reine und angewandte Mathematik*, pages 335–351, 1827. http://www.digizeitschriften.de/dms/img/?PID=PPN243919689_0002%7Clog41.

[19] J. P. De Gua. Sur le nombre des racines Réelles ou Imaginaires, Réelles positives ou Réelles négatives, qui se trouvent dans les équations de tous les degrés. In *Histoire de l'Academie royale des sciences (sec. mémoires)*, pages 72–95. Imprimerie Royale, Paris, 1741.

[20] J. P. De Gua. Recherches du nombre des Racines Réelles ou imaginaires, Réelles positives ou Réelles négatives, qui peuvent se trouver dans les équations de tous les degrés. In *Histoire de l'Academie royale des sciences (sec. mémoires)*, pages 435–494. Imprimerie Royale, Paris, 1741.

[21] B. S. Hourya. Deux moments dans l'histoire du Théorème d'algèbre de Ch. F. sturm. *Revue d'histoire des sciences*, 41(2):99–132, Mar. 1988.

[22] A. Hurwitz. Über den Satz von Budan-Fourier. *Mathematische Annalen*, 71(4):584–591, 1912.

[23] Vilmos Komornik. Another Short Proof of Descartes's Rule of Signs. *The American Mathematical Monthly*, 113(9):829–830, 2006.

[24] P. V. Krishnaiah. A Simple Proof of Descartes's rule of signs. *Mathematics Magazine*, 36(3):135, May–Jun. 1963.

[25] J. L. Lagrange. *Traité de la résolution des équations numériques de tous les degrés*, pages 150–166. Bachelier, Paris, second edition, 1826.

[26] M. N. Obrechkoff. Sur un problème de Laguerre. In *Comptes rendus de l'Académie des Sciences*, pages 102–104. 1923.

[27] F. U. T Æpinus. Démonstration du théorème de Harriot avec une méthode de chercher si une équation algébrique a toutes les racines possibles ou non. In *Histoire de l'Académie Royale des Sciences et des Belles-Lettres de Berlin*, pages 354–366. Haude & Spener, Berlin, 1758.

[28] J. A. Segner. Dissertatio epistolica, qua regulam harriotti... University of Jena, 1728.

[29] J. A. Segner. Démonstration de la règle de Descartes, pour connoitre le nombre des racines affirmatives et négatives qui peuvent se trouver dans les équations. In *Histoire de l'Académie Royale des Sciences et des Belles-Lettres de Berlin*, pages 292–299. Haude & Spener, Berlin, 1756.

[30] B. Szénássy. *History of Mathematics in Hungary until the 20th Century*. Springer-Verlag, Berlin, 1992.

[31] Trélis. Mathématiques. In *Notice des Traveaux de l'Académie du Gard*, pages 195–209. Nisme, 1809.

[32] X. Wang. A Simple Proof of Descartes's Rule of Signs. *The American Mathematical Monthly*, 111(6):525–526, Jun.–Jul. 2004.